\documentclass{amsart}

\usepackage{amsmath,amstext,amsfonts}
\usepackage{amssymb}
\usepackage{latexsym}
\usepackage{epsf}

\newcommand{\Scal}{{\mathcal S}}
\newcommand{\Ucal}{{\mathcal U}}
\newcommand{\B}{{\mathcal B}}

\newcommand{\Z}{{\mathbb Z}}

\newcommand{\Zp}{{\mathbb Z}/p{\mathbb Z}}
\newcommand{\A}{{\mathcal A}}

\newtheorem{Theorem}{Theorem}
\newtheorem{Lemma}{Lemma}
\newtheorem{Conjecture}{Conjecture}
\newtheorem{Corollary}{Corollary}
\newtheorem{Remark}{Remark}

\begin{document}

\title{Sums of dilates in groups of prime order}

\author{Alain Plagne}

\email{plagne@math.polytechnique.fr}
\address{Centre de Math\'ematiques Laurent Schwartz 
\'Ecole polytechnique
91128 Palaiseau Cedex
France}

\maketitle

\bigskip
\begin{center}
{\em 
\`A la m\'emoire du math\'ematicien du d\'esert, Yahya ould Hamidoune,\\
en t\'emoignage d'admiration}
\bigskip 
\end{center}

\begin{abstract}
We obtain a first non-trivial estimate for the sum of dilates problem 
in the case of groups of prime order, by showing that if $t$ is 
an integer different from $0, 1$ or $-1$ and if 
$\A \subset \Zp$ is not too large (with respect to $p$), 
then $|\A+t\cdot \A|>(2+ \vartheta_t)|\A|-w(t)$ for 
some  constant $w(t)$ depending only on $t$ and 
for some explicit real number $\vartheta_t >0$ 
(except in the case $|t|=3$). In the important case $|t|=2$, 
we may for instance take $\vartheta_2=0.08$.
\end{abstract}

\bigskip

\section{Introduction}

It is well known that for two finite sets of integers $\A$ 
and $\B$, their sumset has a cardinality at least
\begin{equation}
\label{basic}
|\A + \B | \geq | \A | + | \B | -1
\end{equation}
and that, except if $\A$ or $\B$ has only one element, 
the equality case in \eqref{basic} implies that $\A$ and $\B$ 
are arithmetic progressions with the same difference; 
in particular, if $|\A|=|\B|$, this implies 
that $\A=\B$ (up to a translation). Freiman's theory \cite{FF} 
then explains that the farther we are from such a situation 
($\A$ and $\B$ close to a given arithmetic progression), 
the bigger the sumset ($\A+\B$) has to be.

The problem of investigating sumsets of dilates takes place in this 
context. We are considering a simplified (one variable) problem 
where we are given a set $\A$ and where, for $\B$, we take a dilation 
of $\A$. In some sense such a set $\B$ is close to $\A$ since 
it is constructed from $\A$, but it is also not that close since 
a dilation for instance does not keep the difference of an arithmetic 
progression unchanged. The problem is then to understand what is going 
on precisely in this situation, by which we mean to measure how far 
from the trivial lower bound (that is, $2|\A |-1$) we are in \eqref{basic}. 

Looking for such an improved bound, in the framework of sets of integers, 
goes back at least to the beginning of the 2000's. In \cite{HP}, 
when dealing with the $3k-3$ theorem of Freiman, the simple estimate 
\begin{equation}
\label{triv}
|\A+t \cdot \A| \geq 3|\A|-2,
\end{equation}
for any integer $t \neq -1,0,1$, was already (and incidentally) proved. 
Then, Nathanson \cite{N} refined this estimate to
$$
|\A + t \cdot \A| \geq \left\lceil \frac72 |\A|- \frac52 \right\rceil,
$$
as soon as $|t| \geq 3$. The case $t=3$ was finally completely solved 
\cite{BB,CSV} (the authors of \cite{CSV} even describe all the cases 
where equality holds)~:
$$
|\A+3 \cdot \A| \geq 4|\A|-4.
$$
In \cite{BB}, Bukh generalizes this result to an arbitrary number of 
summands and proves the nice general lower bound
\begin{equation}
\label{bukh}
|t_1  \cdot \A + \cdots + t_k \cdot \A| \geq (|t_1| + \cdots + | t_k | ) 
|\A | -  o(|\A |)
\end{equation}
as soon as gcd$(t_1, \dots, t_k )=1$, an assumption which is not
restrictive since an affine transformation of $\A$ (by dilations and 
translations) does not change the cardinality of the sumset. 

Returning to the case of two summands which is probably the only one 
where it is reasonable to investigate for precise bounds, in \cite{CHS}, 
the authors obtained
$$
|\A+t \cdot \A| \geq (1+t) |\A| - \left\lceil \frac{t(t+2)}{4} \right\rceil,
$$
if $t$ is prime and $|\A|$ large enough compared to $t$. This was extended 
to the case where $t$ is either a prime power or the product of two primes 
in \cite{DCS}. We shall use these results under the form of a universal 
lower bound (that is, valid for any $\A \neq \emptyset$)
\begin{equation}
\label{pastriv2}
|\A+t \cdot \A| \geq (1+t) |\A| - w(t)
\end{equation}
for some constant $w(t)$. We add that we may take $w(2)=2, w(3)=4$ 
and $w(4)=10$ for instance (using \eqref{triv} and \cite{CHS, DCS}).

The case where the first coefficient is equal to $2$ attracted also some 
energy and the authors of \cite{HR} could prove that, if $t$ is an odd prime,
$$
|2 \cdot \A+ t \cdot \A| \geq (2+t) |\A| - t^2-t+2,
$$
if $|\A|$ is large enough compared to $t$. This was very recently extended 
to the case when $t$ is a prime-power in \cite{ZL}.

In this paper we investigate dilates in $\Z /p\Z$ ($p$ a prime) 
-- a subject which, seemingly, was not yet investigated -- 
and show that a similar phenomenon happens.

\section{The case of $\Z /p\Z$}

While inequality \eqref{basic} in the integers is immediate, 
the Cauchy-Davenport theorem \cite{C} itself is more sophisticated.
In this note, similarly, we are interested in a counterpart to the 
results for integers that were mentioned in the 
Introduction, in the framework of cyclic groups of prime order. 
If one is optimistic, one can hope for a result similar to Bukh's lower 
bound \eqref{bukh}, even with the $o(|\A |)$ term replaced by a constant
depending only on $t$. For the sake of simplicity, and to avoid too many 
technicalities, we shall restrict ourselves to the case of two summands.

\begin{Conjecture}
\label{conj1}
Let $p$ be a prime and $t$ be an integer 
different from 0, 1 or $-1$. Then there exists a constant 
$c(t)$ such that for any set of integers $\A$, the following 
estimate holds
$$
|\A+ t \cdot \A| \geq \min ( (t+1)|\A| - c(t), p).
$$
\end{Conjecture}

Clearly the restriction that only one of the dilation coefficients 
is not equal to $1$ is not restrictive since the cardinality of the sumset 
is invariant by invertible dilation in $\Z /p\Z$ so that 
$| \alpha \cdot \A + \beta \cdot \A| = | \A+ (\alpha^{-1} \beta) \cdot \A |$
when $\alpha$ is non-zero modulo $p$. In the same way, it would not 
be restrictive to impose $|t|<p/2$ (in fact, in this form, 
the conjecture is empty if $t$ is not uniformly bounded with respect 
to $p$).

Notice that, if such a conjecture is true, it implies the result 
in the integers by a cyclification argument (from an additive point 
of view, any set of integers can be seen as a set of elements of 
a cyclic group of large enough prime order). 
Therefore the proof of such a result has to contain the proof 
in the case of integers.

In this paper, we shall not be able to prove Conjecture \ref{conj1} 
but we will rather make a step towards it. To state our result we need to 
introduce a family of auxiliary functions $f_t$ defined as follows 
($t$ is an integral parameter satisfying $t\geq 2$). 

In the case $t=2$, we first introduce 
$$
c_{2}^{(0)}=\frac{1 - 2^{3/2}/3}{2}=0.028595\dots
$$
and then define, for any $0 \leq c \leq 1$, the function $f_2$ as
$$
f_2(c) = \left\{
\begin{array}{l}
\text{the unique solution $x \geq 2$ to the equation $3(1 - cx) = x^{3/2}$,}\\
\hspace{9cm} \text{if $0 \leq c \leq c_2^{(0)}$}, \\
2,\hspace{.5cm} \text{otherwise}.
\end{array}
\right.
$$
This is a decreasing function which assumes the extremal values 
$f_2(0)=3^{2/3}=2.080083\dots$ 
and $f_2(c^{(0)}_2)=2$.

In the case $t \geq 3$, we define
$$
c_t^{(0)}  
= \frac12 \left(1- \frac{2^{3/2}}{(|t|+1) \sin \left(
      \frac{\pi}{|t|+1} \right) } \right).
$$
Notice that $c_t^{(0)}>0$ (except for $t=3$ where $c_3^{(0)}=0$) 
and then define, for any $0 \leq c \leq 1$, the function $f_t$ as
$$
f_{t}(c) = \left\{
\begin{array}{l}
\text{the unique solution $x \geq 2$ to the equation} \\
\hspace{2cm}(|t|+1) \sin \left( \frac{\pi}{|t|+1} \right)(1 - cx) 
= x^{3/2} \sin(\frac{\pi}{x}),\hspace{.5cm} 
\text{if $0 \leq c \leq c_{t}^{(0)}$}, \\
2,\hspace{.5cm} \text{otherwise}.
\end{array}
\right.
$$
All the $f_t$'s functions are again decreasing functions.

We underline the fact that when $t$ tends to infinity, the quantity 
$c_t^{(0)}$ tends towards $(1/2- \sqrt{2} /\pi)=0.04984\dots$ while 
$f_t(0)$ tends towards the unique solution of
$$
x^{3/2} \sin \left( \frac{\pi}{x} \right) = \pi,
$$
that is $2.15409\dots$ These two numerical values are of a certain 
significance explained below.

We are now ready to formulate our main result, which is much more modest 
than Conjecture \ref{conj1} but at least supports this statement
(except for $|t|=3$).

\begin{Theorem}
\label{Th1}
Let $p$ be a prime and $\A$ be a non-empty subset of $\Zp$. 
If $t$ is a prime-power or a product of two primes, 
different from $0, 1$ or $-1$, then
$$
| \A+ t \cdot \A | \geq \min \left( f_{|t|} (c) |\A | - w(t),p \right),
$$
where $c=|\A |/p$, and $w(t)$ the constant depending only on $t$ defined 
in \eqref{pastriv2}.
\end{Theorem}

Notice that if we do not impose to $t$ to be a prime-power or 
a product of two primes, then we can still obtain a non-trivial
result of the form obtained by Bukh.

\begin{Theorem}
Let $\varepsilon >0$. There is an integer $p_0$ depending only on 
$\varepsilon$ such that if $p \geq p_0$ is a prime, $\A$ 
a non-empty subset of $\Zp$ and $t$ is an integer different from 
$0, 1$ or $-1$, then
$$
| \A+ t \cdot \A | \geq \min ( (f_{|t|} (c) - \varepsilon) |\A |, p),
$$
where $c=|\A |/p$.
\end{Theorem}

The proof is the same, mutatis mutandis (that is, essentially 
replacing \eqref{pastriv2} by \eqref{bukh}), as the one of 
Theorem \ref{Th1}.

As an important special case of Theorem \ref{Th1}, we obtain for example 
the following corollary.

\begin{Corollary}
\label{corol}
Let $p$ be a prime and $\A$ be a non-empty subset of $\Zp$ such that 
$|\A | <p/35000$. Let $\sigma$ be equal to either $1$ or $-1$. Then, 
we have
$$
| \A+ 2\sigma \cdot \A | \geq  (2+ \vartheta_2) |\A | - 2,
$$
with $\vartheta_2 =0.08$.
\end{Corollary}

Notice that an equivalent form of the corollary is
$$
| \A+ 2 \cdot \A | \geq  \min ( (2+ \vartheta_2) |\A | - 2, p/17500 -1)
$$

In view of our earlier computation, for large values of $t$ we
would get a $\vartheta_t \sim 0.15409\dots$ as $t$ tends to infinity.
It is therefore a limitation of the method that we cannot hope 
any lower bound coefficient better than $2.16$, say.

This result is clearly demanding for improvements. It should be
possible to extend it to $|t|=3$, with a larger value for $\vartheta_t$ 
and larger sets $\A$, that is, with a larger ratio $|\A | /p$.

The proof relies on a so-called rectification argument, first 
used in Freiman's original proof of his famous theorem \cite{FF}, 
and now standard. It also uses a result of Lev \cite{Lev} and 
its improvement \cite{Lev2}. This proof is presented in the 
following three sections: in Section \ref{secP1} we start with 
the exponential sums argument, then we proceed with the application of 
Lev's results in Section \ref{secP2} and finally, in Section \ref{secP3}, 
we conclude using the rectification trick.

\section{Proof of Theorem \ref{Th1}: the exponential sums method}
\label{secP1}

Let $\A \subset \Zp$ and define $c=|\A|/ p$. 
We define $\Scal = \A+ t \cdot \A$ and define the real number $x$ 
by writing
$$
|\Scal | = x |\A| - w(t).
$$

If $c \geq c_t^{(0)}$ then we simply apply the Cauchy-Davenport theorem 
and get
$$
| \Scal | \geq \min (2| \A | -1, p) \geq \min ( f_t (c) | \A | -w(t), p)
$$
and we are done. From now on and until the end of the proof, 
assume $c \leq c_t^{(0)}$. In this case, we have $f_t (c)|\A|-w(t) <p$ 
so that we are led to prove $|\Scal | \geq f_t(c)|\A | -w(t)$.

If $1_\Ucal$ denotes the characteristic function of a subset $\Ucal$ of
$\Zp$, we write
$$
\widehat{1_\Ucal} (x) = \sum_{u\in \Ucal} \exp (2 \pi i ux/p),
$$
for the discrete Fourier transform of $1_\Ucal$ and use 
the exponential sums counting method, which yields
$$
p|\A|^2 = 
\sum_{r=0}^{p-1} \widehat{1_\A} (r) \widehat{1_{t \cdot \A}} (r) 
\overline{ \widehat{1_\Scal} (r)}.
$$
Then, we proceed with the classical calculations 
(using Cauchy-Schwartz inequality, Parseval identity,\dots):
\begin{eqnarray*}
p|\A|^2 &  =  & |\A |^2 | \Scal| + 
\sum_{r=1}^{p-1} \widehat{1_\A} (r) \widehat{1_{t \cdot \A}} (r) 
\overline{ \widehat{1_\Scal} (r)} \\
	& \leq & |\A|^2|\Scal| + \max_{1\leq r \leq p-1} 
| \widehat{1_\A} (r)|  
\sum_{r =1}^{p-1} |\widehat{1_{t \cdot \A}} (r)| |\widehat{1_\Scal} (r)| \\
       & \leq & |\A|^2|\Scal| + \max_{1 \leq r \leq p-1} 
| \widehat{1_\A} (r)| \ 
\left( \sum_{r =1}^{p-1} |\widehat{1_{t \cdot \A}} (r)|^2 \right)^{1/2} 
\left( \sum_{r =1}^{p-1} |\widehat{1_\Scal} (r)|^2 \right)^{1/2} \\
       & \leq & |\A|^2|\Scal| + \max_{1 \leq r \leq p-1} | 
\widehat{1_\A} (r)|\  
\left( \sum_{r=0}^{p-1} |\widehat{1_{\A}} (tr)|^2 \right)^{1/2} 
\left( \sum_{r=0}^{p-1} |\widehat{1_\Scal} (r)|^2 \right)^{1/2} \\
 & \leq & |\A|^2|\Scal| + \max_{1 \leq r \leq p-1} | \widehat{1_\A} (r)|\ 
\sqrt{p|\A|} \sqrt{p|\Scal|} \\
 &  =   & x |\A|^3 +   p  \sqrt{|\A| |\Scal |}\ \max_{1 \leq r \leq p-1} 
| \widehat{1_\A} (r)|
\end{eqnarray*}
from which it follows the following lower bound for the Fourier bias
of $\A$,
$$
\max_{1 \leq r \leq p-1} | \widehat{1_\A} (r)| \geq 
\frac{p - x|\A|}{ p \sqrt{x}} |\A| = \frac{1 - xc}{\sqrt{x}} |\A|. 
$$

Again, since the cardinality of a set is independent of the invertible 
dilation in $\Zp$, we may assume without loss of generality that 
this maximum is attained for $r=1$.

\section{Proof of Theorem \ref{Th1} continued: a lemma by Lev}
\label{secP2}

We now use (a special case of) Theorem 1 of \cite{Lev} and Corollary 1 
of \cite{Lev2} which we state here as a single unified lemma for our 
purpose. We shall make use of the function defined on $(0,\pi]$
$$
g(u)= \frac{\sin u}{u}
$$
and of its reciprocal function $g^{-1}$ which is decreasing on $[0,1)$.

\begin{Lemma}
With the preceding notation, define 
$$
\eta = \frac{\left|\widehat{1_\A}(1)\right|}{|\A|}.
$$ 
Then, for any $\beta \in (0,1/2]$, there exists an interval $I$ 
of length $\beta p$ in $\Zp$ such that
$$
|\A \cap I | \geq M (\beta, \eta) | \A|
$$
where
$$
M (\beta, \eta) = \max \left(
\frac{\eta +1 - 2 \cos \pi \beta}{2(1- \cos \pi  \beta)}, 
\frac{\pi \beta }{g^{-1}(\eta g(\pi \beta ))} \right).
$$
\end{Lemma}

In view of what will be needed later in the proof, we now state
an important remark.

\begin{Remark}
Assume that $0 \leq \eta \leq 1/\sqrt{2}$ and $0 \leq \beta \leq 1/3$.
If
$$
\beta^{-1} \left(\frac{\eta +1 - 2 \cos \pi \beta}{2(1- \cos \pi
    \beta)}\right) \geq 2
$$
then $\beta \geq 1/4$.
If 
$$
\frac{\pi }{g^{-1}(\eta\ g(\pi \beta ))} \geq 2
$$
then $\beta \leq 1/4$.
\end{Remark}

Indeed the first inequation leads to
$$
\cos \pi \beta \leq \frac{1+\eta-4\beta}{2(1-2\beta)} 
\leq \frac{1+1/\sqrt{2}-4 \beta}{2(1-2\beta)},
$$
which is easily shown to imply $\beta \geq 1/4$.

For the second inequation, applying $g$, it implies
$$
\eta g (\pi \beta) \geq g \left( \frac{\pi}{2} \right)= \frac{2}{\pi},
$$
and therefore (again due to $\eta \leq 1/\sqrt{2}$) that  
$$
g (\pi \beta) \geq \frac{2 \sqrt{2}}{\pi}= g \left( \frac{\pi}{4} \right)
$$
and we obtain $\beta \leq 1/4$.

\section{Proof of Theorem \ref{Th1} concluded: the rectification}
\label{secP3}

We apply Lev's lemma with $\beta=1/(|t|+1)$ and the value we just obtained 
for $\eta$, namely $\eta= (1 - xc)/\sqrt{x}$. We obtain an interval $I$ 
of length $[p/(|t|+1)] \leq (p-1)/(|t|+1)$ 
such that $\A_0= \A \cap I$ contains at least
\begin{equation}
\label{wok}
| \A_0 | \geq B_t (x,c) | \A |
\end{equation}
elements, where we denote by $B_t(x,c)$ the lower bound 
$$
B_t(x,c) = M (\beta, \eta)= M \left( \frac{1}{|t|+1},
\frac{1-xc}{\sqrt{x}} \right).
$$

The point is now to notice that the sumset $\A_0 + t \cdot \A_0$
is Freiman isomorphic to the same set seen in the integers. Indeed all 
the elements in the sumset belong to an interval of integers of length 
$(1+|t|) [p/(|t|+1)] \leq p-1 <p$, and therefore two sums are equal if 
and only if they are equal modulo $p$. Thus, we can apply to the sumset 
$\A_0$ the lower bound derived in the case of integers \eqref{pastriv2} 
and get, in view of \eqref{wok},
\begin{eqnarray*}
x |\A |- w(t) =|\Scal | =|\A+ t \cdot \A | & \geq & | \A_0 + t \cdot \A_0| \\
		  & \geq &  (|t|+1) |\A_0| - w(t) \\
		  & \geq &  (|t|+1) B_t (x,c) |\A | -w(t).
\end{eqnarray*}
We finally obtain
\begin{equation}
\label{borninf}
x \geq (|t|+1) B_t (x,c). 
\end{equation}
In order to know which bound is worth here (for the maximum), 
we use the remark stated just after Lev's lemma.

In the case $|t|=2$, we use the first argument in the maximum defining $M$. 
It then follows from \eqref{borninf}
$$
x \geq  3B_2(x,c) = 3 \left( \frac{1 - xc}{\sqrt{x}} \right)
$$
or equivalently $x \geq f_2(c)$ by definition of the function $f_2$ and the 
result is proved.

As for proving Corollary \ref{corol}, we only have to solve 
$f_2(c)=2.08$ to get
$$
c = \frac{1-2.08^{3/2}/3}{2.08}=0.0000209607\dots =
\frac{1}{34410.7\dots}> \frac{1}{35000}.
$$

In the case $|t| \geq 3$, we use the second argument in the maximum 
defining $M$, which gives similarly
$$
x \geq \frac{\pi}{g^{-1}(\eta\ g(\pi/(|t|+1)) )}. 
$$
Applying $g$, this can be rewritten as
$$
\eta\ g \left( \frac{\pi}{|t|+1} \right) \leq g \left( \frac{\pi}{x} \right),
$$
that is
$$
(|t|+1) \sin \left( \frac{\pi}{|t|+1} \right) (1-xc) \leq 
x^{3/2} \sin \left( \frac{\pi}{x} \right).
$$
And by definition, this implies
$$
x \geq f_t (c)
$$
and the result follows.

\bigskip
\bigskip

{\bf Acknowledgements:} This note was started in Paris in January 2011 
after a long discussion with Yahya ould Hamidoune on the case of integers. 
The author had just the time to announce him the (quantitative version 
of the) result before he stopped any mathematical activity. 
The article was then mainly written down while the author was enjoying 
an invitation to the Isaac Newton Institute for Mathematical Sciences 
in Cambridge at the occasion of the {\em Discrete Analysis} Programme, 
February 2011. He thanks the organizers for this kind invitation and 
the good working conditions they offered.

\bigskip
\bigskip


\bigskip
\bigskip




\begin{thebibliography}{10}


\bibitem{BB} B. Bukh, {\it Sums of dilates}, 
Combin. Probab. Comput. 17 (2008), no. 5, 627--639.

\bibitem{C} A.-L. Cauchy, {\it Recherches sur les nombres}, 
J. \'Ecole polytechnique 9 (1813), 99--123.

\bibitem{CHS} J. Cilleruelo, Y. ould Hamidoune, O. Serra, 
{\it On sums of dilates}, 
Combin. Probab. Comput.  18  (2009), 871--880.

\bibitem{CSV} J. Cilleruelo, M. Silva, C. Vinuesa, 
{\it A sumset problem}, J. Combin. Number Th. 2 (2010).

\bibitem{DCS} S.-S. Du, H.-Q. Cao, Z.-W. Sun, 
{\it On a sumset problem for the integers}, arXiv:1011.5438.

\bibitem{FF} G. A. Freiman, {\it Foundations of a structural
theory of set addition}, Trans. AMS Monographs {\bf 37}, AMS, 1973.

\bibitem{HP} Y. ould Hamidoune, A. Plagne, 
{\it A generalization of Freiman's $3k-3$ theorem}, 
Acta Arith. 103 (2002), no. 2, 147--156.

\bibitem{HR}  Y. ould Hamidoune, J. Ru\'e, {\it On dilates sums}, 
Combin. Probab. Comput. 20 (2011), 249--256.

\bibitem{Lev} V.F. Lev, {\it Distribution of points on arcs}, Integers
  5 (2005), 8pp.

\bibitem{Lev2} V.F. Lev, {\it More on points and arcs}, Combinatorial 
number theory,  347--350, de Gruyter, Berlin, 2007.

\bibitem{ZL} Z. Ljujic, {\it A lower bound for the size of a sum of dilates},
arXiv 11.01.5425.

\bibitem{N} M. Nathanson, {\it Inverse problems for linear forms over
    finite sets of integers}, J. Ramanujan Math. Soc. 23 (2008),
  no. 2, 151--165.

\end{thebibliography}
\end{document}